\documentclass[12pt]{article}
\usepackage{amssymb,amsmath,amsfonts}

\begin{document}

\pagestyle{plain}

\title{ On $p$-separability of subgroups of free metabelian groups}

\author{Valerij G. Bardakov}


\maketitle

\begin{abstract}
We prove that every free metabelian non--cyclic group
    has a finitely generated isolated subgroup
which is not separable in the class of nilpotent
groups.

As a corollary we prove that for every prime number $p$
an arbitrary free metabelian non--cyclic group has a
finitely generated $p'$--isolated subgroup which is
not $p$--separable.

\noindent
{\it Mathematics Subject Classification:} 20F16, 20F14, 20F10\\
\noindent
{\it Key words and phrases:} free metabelian group, nilpotent
group, isolated subgroup, $p$--separable subgroup, Magnus
representation.
\end{abstract}


Following A. I. Mal'cev [1] we say that a subgroup $H$ of a group
$G$ is {\it separable in a class of groups} $\mathcal{K},$ if for
each $g \in G \setminus H$ there is a homomorphism $\varphi$ of
$G$ to some  group from $\mathcal{K}$ such that $\varphi(g) \not
\in \varphi(H)$.  If $\mathcal{K}$ is the class of all finite
groups (resp. finite $p$--groups), the corresponding notion of
$\mathcal K$-separability is called the {\it finite separability}
(resp. {\it finite $p$--separability}). The problem of finite
separability is closely related to the generalized word problem
[1]. (The {\it generalized word problem} for $H$ in $G$ asks for
an algorithm that decides whether or not the elements of $G$ lie
in $H$.)

Let $p$ be a prime number. Recall that a subgroup $H$ of $G$ is
called $p'$--{\it isolated}, if for every prime $q \not =p$ and
for every $g \in G$ the condition $g^q \in H$ implies that $g \in
H$. E.~D.~Loginova [2, \S~3] proved that in each finitely
generated nilpotent group every $p'$--isolated subgroup is
$p$--separable.

D.~I.~Moldavanskii suggested that the latter fact is true
in every free group.

{\bf Problem} ([3, Problem 15.60]). {\it  Is it true
that any finitely generated $p'$--iso\-la\-ted
subgroup of a free group is separable in the class of
finite $p$--groups?  It is easy to see that this is
true for cyclic subgroups.}

The paper [4] gives the negative answer to Moldavanskii's problem.
More precisely, it is proved there that in each non--abelian free group
$F$ there is an isolated (and therefore a
$p'$--isolated) finitely generated subgroup which is
not $p$--separable. In particular, the result by
E.~D.~Loginova cannot be generalized to the class of
absolute free groups.

In the present article we shall demonstrate that the
result by E.~D.~Loginova cannot be generalized to the
class of soluble groups, too. We prove the following
theorem.

{\bf Theorem.} {\it Every free metabelian
non--cyclic group contains a finitely generated
isolated subgroup which is not separable in the class
of nilpotent groups.}

Since every  finite  $p$--group is nilpotent [5, p.
162], and every isolated subgroup is $p'$--isolated
for each prime $p$, our theorem implies

{\bf Corollary.} {\it Every  free metabelian
non--cyclic group contains a finitely generated
$p'$--isolated subgroup which is not $p$--separable.}

\vskip 20pt

 \centerline{\S~1. {\bf The Magnus representation of free metabelian groups}}

\vskip 12pt

Let $F$ denote a free group of rank $n$ and $F''$
the second commutator subgroup of $F$. Then $\Phi = F/F''$
is a free metabelian group of rank $n$.
W.~Magnus [6, Chapter I, \S~4] constructed a
faithful $2\times 2$ matrix representation of $\Phi$.
S.~Bachmuth [7] studied the properties of the Magnus
representation. For the reader's convenience
we reproduce some properties of the Magnus representation.

Let $\Phi$ be freely generated by $x_1, x_2,\ldots, x_n$. Let
$s_1, s_2,\ldots, s_n$, $t_1, t_2,\ldots, t_n$ be commuting
indeterminates. Then $\Phi$ has a faithful representation by
$2\times 2$ matrices defined via the correspondence
$$
\mu : x_i \longrightarrow
\begin{pmatrix}
s_i &  t_i    \\
0 & 1  \\
\end{pmatrix}
,~~~
i = 1, 2,\ldots, n.
$$

Let $c = s_1^{j_1} s_2^{j_2}\ldots s_n^{j_n}$, where $j_1, j_2,\ldots, j_n$
are arbitrary integers. Suppose that
$$
\gamma_1, \gamma_2,\ldots, \gamma_n
\in \mathbb{Z}[s_1^{\pm 1}, s_2^{\pm 1},\ldots, s_n^{\pm 1}].
$$
Then the matrix
$$
\begin{pmatrix}
c &  \sum\limits_{i=1}^n \gamma_i t_i    \\
0 & 1  \\
\end{pmatrix}
$$
is in the image of $\mu $, or, equivalently,
it is a product of matrices
$$
\begin{pmatrix}
s_i &  t_i    \\
0 & 1  \\
\end{pmatrix}
^{\pm 1}
$$
provided that the $\gamma_i$ satisfy the identity
$$
\sum\limits_{j=1}^n \gamma_i (1 - s_i) = 1 - c
$$
and vice versa (see [7, Lemma 1]).

\vskip 20pt

 \centerline{\S~2. {\bf Proof of the Theorem}}

\vskip 12pt

Let $\Phi$ be a free metabelian group with free generators
$$
 x,~~y,~~ z_i~~(i \in I),
$$
where $I$ is some index set (possibly empty).
Let $H$ be a subgroup of $\Phi$ generated by
the elements
$$
a = x [y, x],~~~ b = y, ~~~z_j~~ (j \in J),
$$
where $[y, x] = y^{-1} x^{-1} y x$ and $J$ is a subset of $I$.

We claim that $H$ is a proper subgroup of $\Phi$ and,
in particular, the element $x$ is not in $H$.
To do this let us consider a map $\tau$ from $\Phi$ to
the symmetric group $S_3$ defined via
$$
\tau : \left\{
\begin{array}{ll}
x \longmapsto (1 2), &  \\
y \longmapsto (2 3), & \\
z_{i} \longmapsto 1 & \mbox{ if }~~ i \in I.
\end{array} \right.
$$
Since $S_3$ is a metabelian group, $\tau $ is a homomorphism.
It is easy to see that $\tau(a) = \tau(b) = (2 3)$, i.~e.,
$\tau (H) = \langle (2 3) \rangle \simeq \mathbb{Z}_2$ and $\tau (x) \not \in \tau (H).$
Hence $x \not \in H,$ as required.

Now let us prove that the subgroup $H$ is not separable in the class of nilpotent groups.
Consider the lower central series of $\Phi$
$$
\Phi = \gamma_1 \Phi \geq \gamma_2 \Phi \geq \gamma_3 \Phi \geq
\ldots,
$$
where
$$
\gamma_{i+1} \Phi = [ \gamma_i \Phi, \Phi ],~~~i = 1, 2,
\ldots ,
$$
and the set of homomorphisms
$$
\varphi_n : \Phi \longrightarrow \Phi/\gamma_n \Phi,~~~ n = 1, 2, \ldots ,
$$
to the nilpotent metabelian groups.
Note that the image of $x$
under all these homomorphism is in the image of $H$.
Indeed, it is easy to check that the subgroup
generated by the elements
$$
\varphi_n(a), ~~\varphi_n(b), ~~\varphi_n(z_j), ~~ j\in J,
$$
and the subgroup generated by the elements
$$
\varphi_n(x), ~~\varphi_n(y), ~~\varphi_n(z_j), ~~j \in J,
$$
are equal modulo the commutator
subgroup $\Phi' = \gamma_2 \Phi.$

Then (see [8, Theorem 31.2.5]) the group
$$
\langle \varphi_n(a), ~~\varphi_n(b),~~ \varphi_n(z_j)~~ (j \in J) \rangle
$$
is equal to the group
$$
\langle \varphi_n(x),~~ \varphi_n(y),~~ \varphi_n(z_j)~~ (j \in J) \rangle.
$$
Therefore $\varphi_n(x) \in \varphi_n(H)$ for all natural $n$.
It then follows that  $H$ is not separable from element
 $x$ in the class of free nilpotent groups and hence in the class
 of all nilpotent groups.

To complete the proof of the Theorem we must prove that
$H$ is isolated in $\Phi$.
For the sake of simplicity,
we consider the case when
the index set $I$ is empty, i.~e., $\Phi$ is
generated by $x$ and $y$.
The proof  in the general
case is similar.

Now consider the matrices
$$
X = \mu(x) =
\begin{pmatrix}
\alpha &  e    \\
0 & 1  \\
\end{pmatrix},~~~
Y = \mu(y) =
\begin{pmatrix}
\beta &  f    \\
0 & 1  \\
\end{pmatrix},
$$
where $\{e,f\}$ is a base of a free module of rank $2$
over the ring $R = \mathbb{Z}[\alpha^{\pm 1},
\beta^{\pm 1}]$. Lemma 1 from [7] implies that the matrix
$$
\begin{pmatrix}
\alpha^{i_1} \beta^{i_2} &  \gamma_1 e + \gamma_2 f    \\
0 & 1  \\
\end{pmatrix},
$$
where $i_1, i_2 \in \mathbb{Z},$ $\gamma_1, \gamma_2
\in R$, is in $\mu(\Phi)$ if and only if
$$
\gamma_1 (1 - \alpha ) + \gamma_2 (1 - \beta ) = 1 - \alpha^{i_1}
\beta^{i_2}.
$$
Note that there exist an epimorphism $\psi : \Phi
\longrightarrow H$ that sends
$x$ to $a$ and $y$ to $b$. In this case we
have the following diagram
$$
\begin{array}{ccc}
\Phi & \stackrel{\psi}{\longrightarrow} & H    \\
\mu \downarrow~~ &  & ~~~\downarrow \mu \\
M_2(Q) & \stackrel{\psi^*}{\longrightarrow}  & M_2(Q)    \\
\end{array}
$$
where $Q = R[e, f]$. In order to define the map $\psi^{\ast}$, which makes this diagram to
commutative we will find the images of generators of $H$ under $\mu $. Note that
$$
X^{-1} =
\begin{pmatrix}
\alpha^{-1} &  -\alpha^{-1} e    \\
0 & 1  \\
\end{pmatrix},~~~
Y^{-1} =
\begin{pmatrix}
\beta^{-1} &  -\beta^{-1} f    \\
0 & 1  \\
\end{pmatrix},
$$
which trivially implies that
$$
\mu(a) = X [X, Y] =
\begin{pmatrix}
\alpha & (2 - \beta^{-1}) e + (1 - \alpha )\beta^{-1} f    \\
  &   \\
0 & 1  \\
\end{pmatrix},
$$
$$
\mu(b) = Y.
$$
Write $A$ for $\mu(a)$ and $B$ for $\mu(b).$
We construct the homomorphism $\psi^{\ast}$
via
$$
\psi^{\ast} (X) = A,~~~\psi^{\ast} (Y) = B.
$$
It is clear that the homomorphism
$\lambda $
of free module $M$ defined as follows
$$
\lambda : \left\{
\begin{array}{ll}
e \longmapsto (2 - \beta^{-1}) e + (1 - \alpha )\beta^{-1} f, &  \\
f \longmapsto f, & \\
\end{array} \right.
$$
induces
the homomorphism $\psi^* : M_2(Q) \longrightarrow M_2(Q).$
Then the map $\psi^{\ast}$ takes a matrix
$$
S =
\begin{pmatrix}
\alpha^{i_1} \beta^{i_2} &  \gamma_1 e + \gamma_2 f    \\
0 & 1  \\
\end{pmatrix}
 \in M_2(Q)
$$
to the matrix
$$
\psi^{\ast} (S) =
\begin{pmatrix}
\alpha^{i_1} \beta^{i_2} &
\gamma_1 (2 - \beta^{-1}) e + ((1 - \alpha) \beta^{-1} \gamma_1 + \gamma_2) f    \\
  &   \\
0 & 1  \\
\end{pmatrix}.
$$
Now it is easy to see that this map is a homomorphism.

 The matrix of the endomorphism $\lambda $ in the base $e$, $f$
is
$$
\begin{pmatrix}
2 - \beta^{-1} &  (1 - \alpha ) \beta^{-1}    \\
0 & 1  \\
\end{pmatrix}.
$$
Its determinant is equal to $2 - \beta^{-1} = (2\beta - 1)/\beta
$. Hence $\lambda$ is not automorphism of $M$ being considered as
a module over $\mathbb{Z}[\alpha^{\pm 1}, \beta^{\pm 1}].$ The
module $M$ can be placed into the vector space $\widetilde{M}$
over field of rational fractions $\mathbb{Q}(\alpha, \beta)$ with
the base $e$, $f.$ Now $\lambda $ induces  an automorphism of
$\widetilde{M}$; we shall denote the said induced map by
$\widetilde{\lambda}.$ Hence $\widetilde{\lambda} \in
\text{Aut}(\widetilde M)$ does possess the inverse and
$$
\widetilde{\lambda}^{-1} : \left\{
\begin{array}{ll}
\displaystyle e \longmapsto \frac{\beta}{2\beta - 1} e + \frac{\alpha - 1}{2\beta - 1} f, &  \\
  &   \\
f \longmapsto f. & \\
\end{array} \right.
$$
The map $\widetilde{\lambda}^{-1}$ induces a map $\psi_1^{\ast}$
which sends the matrix
$$
S =
\begin{pmatrix}
\alpha^{i_1} \beta^{i_2} &  \gamma_1 e + \gamma_2 f    \\
0 & 1  \\
\end{pmatrix}
\in M_2(Q)
$$
to the matrix
$$
\psi_1^{\ast} (S) =
\begin{pmatrix}
\alpha^{i_1} \beta^{i_2} &
\widetilde{\lambda}^{-1} (\gamma_1 e + \gamma_2 f)    \\
0 & 1  \\
\end{pmatrix}.
$$
However, it is not true in general that
the latter matrix belongs to the ring
$M_2(Q)$. The following simple
lemma gives a criterion when it does
belong to $M_2(Q).$

{\bf Lemma 1.} {\it The matrix $S \in \mu(\Phi)$ belongs
to the subgroup $\mu(H)$ if and only if $\psi_1^{\ast
}(S) \in \mu(\Phi)$.}

Let us return back to the proof of the Theorem. Consider
the matrix $S$ in $M_2(Q) \cap \mu(\Phi)$. By Bachmuth's lemma above,
the following
equality
$$
\gamma_1 (1 - \alpha ) + \gamma_2 (1 - \beta ) = 1 - c,
$$
where $c = \alpha^{i_1} \beta^{i_2}$ is true.
Assume that $S^m \in \mu(H)$ for some positive integer
$m$.

The following lemma is obvious.

{\bf Lemma 2.} {\it For every natural $m$ and every matrix
$$
S =
\begin{pmatrix}
c &  \gamma_1 e + \gamma_2 f    \\
0 & 1  \\
\end{pmatrix}
 \in M_2(Q),
$$
where $c = \alpha^{i_1} \beta^{i_2}$ for some $i_1, i_2 \in \mathbb{Z}$,
the following formulas are true}
$$
S^m =
\begin{pmatrix}
c^m &  (1 + c + \ldots + c^{m-1}) (\gamma_1 e + \gamma_2 f)    \\
  &   \\
0 & 1  \\
\end{pmatrix},
$$
$$
S^{-m} =
\begin{pmatrix}
c^{-m} &  -(c^{-1} + c^{-2} + \ldots + c^{-m}) (\gamma_1 e + \gamma_2 f)    \\
  &   \\
0 & 1  \\
\end{pmatrix}.
$$

By Lemma 2 matrix $S^m$ is equal to
$$
S^m =
\begin{pmatrix}
c^m &  (1 + c + \ldots + c^{m-1}) (\gamma_1 e + \gamma_2 f)    \\
  &   \\
0 & 1  \\
\end{pmatrix}.
$$
Then by Lemma 1 the matrix
$$
\psi_1^{\ast }(S^m) =
\begin{pmatrix}
\displaystyle c^m &  \displaystyle (1 + c + \ldots + c^{m-1}) \left(\frac{\gamma_1 \beta }{2\beta - 1} e +
\displaystyle \left(\frac{\gamma_1 (\alpha - 1)}{2\beta - 1} + \gamma_2 \right) f\right)    \\
  &   \\
0 & 1  \\
\end{pmatrix}
$$
is in $M_2(Q)$. However, it is possible only in the case when
$(1 + c + \ldots + c^{m-1}) \gamma_1 \beta$ is a multiple of  $2 \beta - 1$.
Note that the polynomials  $\beta $ and $2 \beta - 1$ are relatively
prime
and hence the polynomial
$1 + c + \ldots + c^{m-1}$
is relatively prime with $2 \beta - 1$. Indeed, our
first statement is evident while the second one follows from
the fact that $c^{m-1} = (\alpha^{i_1} \beta^{i_2})^{m-1}$ and
hence cannot be a multiple of $2$.
Therefore, $\gamma_1$ is
a multiple of $2 \beta - 1$; but this means that
$\psi_1^{\ast }(S) \in M_2(Q)$.
We are going to show that  $\psi_1^{\ast }(S) \in
\mu(\Phi)$; then in view of Lemma 1 this will mean that
 $S \in \mu(H)$.

Consider the matrix
$$
\psi_1^{\ast }(S) =
\begin{pmatrix}
c & \displaystyle \frac{\gamma_1 \beta }{2\beta - 1} e +
\left(\frac{\gamma_1 (\alpha - 1)}{2\beta - 1} + \gamma_2 \right)f    \\
  &   \\
0 & 1  \\
\end{pmatrix}.
$$
By Bachmuth's lemma, this matrix is in $\mu(\Phi)$ provided
that
$$
\frac{\gamma_1 \beta }{2\beta - 1} (1 - \alpha ) +
\left( \frac{\gamma_1 (\alpha - 1)}{2\beta - 1} + \gamma_2 \right)
(1 - \beta ) = 1 - c.
$$
It is easy to check that this equality is equivalent to
$$
\gamma_1 (1 - \alpha ) (2\beta - 1) +
\gamma_2 (1 - \beta) (2\beta - 1) = (1 - c) (2 \beta - 1).
$$
This implies that
$$
\gamma_1 (1 - \alpha ) +
\gamma_2 (1 - \beta ) = 1 - c.
$$
Since $S \in \mu(\Phi),$ it follows that the latter equality is true.

The Theorem is proven.

{\bf {\it Acknowledgements.}} I am very grateful to
M.~V.~Neshchadim for very stimulating discussions. I also appreciate advice
from Vladimir Tolstykh concerning
certain technical matters. Special thanks
goes to the participants of the seminar ``Evariste Galois'' at
Novosibirsk State University for their kind attention to my work.

\vskip 20pt


\vskip 24pt

\centerline{REFERENCES} \vskip 12pt
\begin{enumerate}
\item
    A. I. Mal'cev, On homomorphisms onto finite groups, Uchen. Zapiski Ivanovsk. ped. instituta,
18, N~5 (1958), 49--60 (also in ``Selected papers'', Vol. 1, Algebra, 1976, 450--462) (Russian).
\item
E.~D.~Loginova, Residual finiteness of the free product of two
groups with commuting subgroups, Sib. Mat. Zh. 40, N~2 (1999),
395-407 (Russian).
\item
   The Kourovka Notebook (Unsolved problems
in group theory), 15th ed., Institute of Mathematics SO RAN, Novosibirsk, 2002.
\item
    V. G. Bardakov,
On D. I. Moldavanskii's Question About p-Separable Subgroups of a
Free Group,
 Siberian Mathematical Journal, 45, No~3 (2004), 416--419; see also arXiv :
math.GR/0404292.
\item
        M.~I.~Kargapolov, Yu.~I.~Merzljakov, Fundamentals of the Theory of Groups, New York:
Springer, 1979.
\item
R. C.~Lyndon, P. E.~Schupp, Combinatorial group theory,
Ergebnisse der Mathematik und ihrer Grenzgebiete,
Springer-Verlag, Berlin-New York, 1977.
\item
    S. Bachmuth, Automorphisms of free metabelian groups, TAMS,
118, N~6 (1965), 93--104.
\item
    H. Neumann, Varieties of Groups, Ergebnisse der Mathematik und ihrer Grenzgebiete, Band 37,
Springer--Verlag, Berlin--Heidelberg--New York, 1967.

\end{enumerate}

\bigskip
\bigskip
\noindent
Author address:

\bigskip

\noindent
Valerij G. Bardakov\\ Sobolev Institute of Mathematics,\\
Novosibirsk, 630090, Russia\\ {\tt bardakov@math.nsc.ru}

\end{document}